\documentclass[12pt]{amsart}
\usepackage{amsmath}
\usepackage{amsfonts,amssymb,amsthm}

\newtheorem{theorem}{Theorem}

\theoremstyle{remark}
\newtheorem{remark}[theorem]{Remark}

\newcommand{\vol}{{\mathrm{vol}\,}}

\begin{document}


\title[Ellipsoid estimates]{Simple estimates for ellipsoid measures}

\author{Igor Rivin}

\address{Department of Mathematics, Temple University, Philadelphia}

\curraddr{Mathematics Department, Princeton University}

\email{rivin@math.temple.edu}

\thanks{The author is supported by the NSF DMS}

\date{\today}

\keywords{ellipsoid, mean curvature, surface area, John ellipsoid}

\subjclass{52A38; 11P21}

\begin{abstract}
We write down estimates for the surface area, and more generally, integral
mean curvatures of an ellipsoid $E$ in $\mathbb{E}^n$ in terms of the lengths
of the major semi-axes. We give applications to estimating the area of
parallel surfaces and volume of the  tubular neighborhood of $E,$ to the
counting of lattice points contained in $E$ and to estimating the shape of the
John ellipsoid of a convex body $K.$
\end{abstract}

\maketitle

\section*{Introduction}
Consider an ellipsoid $E \subset \mathbb{E}^n.$ Such an ellipsoid is defined
as the set
\begin{equation}
\label{edef}
E = \{x\in \mathbb{E}^n \quad | \quad \|A x\| = 1\},
\end{equation}
where $A$ is a non-singular linear transformation of $\mathbb{E}^n.$
Remark that 
\[\|A x\| = \langle A x, A x\rangle = \langle x, A^t A x\rangle.
\]
The matrix $Q = A^t A$ is a positive definite matrix, with eigenvalues
$\lambda_1, \dots, \lambda_n,$ whose (positive) square roots $\sigma_1, \dots,
\sigma_n$ are the so-called \emph{singular values} of $A.$ Their geometric
significance is that the \emph{major semi-axes} of $E$ are the quantities $a_i
= 1/\sigma_i.$ After an orthogonal change of coordinates, we can write
\begin{equation}
\label{edef2}
E = \{x \in \mathbb{E}^n \quad | \quad \sum_{i=1}^n x_i^2 \lambda_i = 
\sum_{i=1}^n x_i^2 \sigma_i^2 = \sum_{i=1}^n \dfrac{x_i^2}{a_i^2} = 1\}.
\end{equation}
It is evident from  \eqref{edef} and \eqref{edef2} that the volume of $E$ is
given by 
\begin{equation}
\label{evol}
\vol E = \dfrac{\kappa_n}{\det A} = \kappa_n \prod_{i=1}^n a_i,
\end{equation}
where $\kappa_n$ is the volume of the unit ball in $\mathbb{E}^n.$
The formula \eqref{evol} is deceptive, in that, as is quite well known, there
is no simple expression (in terms of the major semi-axis lengths) for the
\emph{surface area} of $E,$ even in $2$ dimensions. In this note we write down
an estimate for the surface area and \emph{integral mean curvatures of $E$} in
terms of the major semi-axes. These estimate differ from the true values by a
factor which depends only on the dimension $n.$ It should be noted that the
estimate on the dimensional factor is extremely crude, but the methods used to
get sharper results are quite different from those used in this note, and so we
postpone that to a different note \cite{rellip2}. In any event, these estimate
also allow us to estimate the volume of the set
\[
E_\rho = \{ x \quad | \quad d(x, E) \leq \rho\},
\]
as well as of the surface area of $\partial E_\rho.$ Despite appearances,
$E_\rho$ is usually \emph{not} an ellipsoid.

Our estimates for the integral mean curvatures (which are defined in the next
section of this note) will have the following form:
\begin{equation}
\label{meanest}
c_{n, i} \leq \dfrac{\mathcal{M}_i(\partial E)}{s_{n-1-i}(a_1, \dots, a_n)} \leq C_{n,
  i},
\end{equation}
where $C_{n, i}/c_{n, i} \leq n^{(n-i)/2},$
and $s_k$ is the $k$-th elementary symmetric function, defined by
\[
\prod_{i=1}^n (x + a_i) = \sum_{k=0}^n x^k s_{n-k}(a_1, \dots, a_n).
\] It should be noted that $\mathcal{M}_0$ is the surface area of $\partial
E.$ 
The estimates for $\vol E_\rho$ and $\vol_{n-1} \partial E_\rho$ are given in
Theorem \ref{parest}.

\section{Integral geometry}
The reference for the results recalled in this section is Santalo's book
\cite{santalo}. Let $K$ be a convex body in $\mathbb{E}^n;$ we will initially
assume that $\partial K$ is of smoothness  at least $C^2.$ The \emph{$k$-th
  mean curvature $m_k(\partial K, x)$} is defined as 
\[
\binom{n-1}{k} m_k(x) = s_k(k_1(x), \dots, k_{n-1}(x)),
\]
where $k_1, \dots, k_{n-1}$ are the \emph{principal curvatures} of $\partial
K$ at the point $x.$
The \emph{$k$-th integral mean curvature $\mathcal{M}_k$} is defined as 
\[
\mathcal{M}_k = \int_{\partial K} m_k(x) d A,
\]
where $d A$ is the usual surface measure on $\partial K.$ In particular, since
$m_0 = 1,$ it follows that $\mathcal{M_0} = \vol_{n-1} \partial K.$

The properties of $\mathcal{M}_k$ we will use are:
\begin{enumerate}
\item
The area of $\partial K_\rho$ is given by (\cite[(13.43)]{santalo}:
\begin{equation}
\label{arho}
\vol_{n-1} \partial K_\rho = \sum_{k=0}^{n-1} \binom{n-1}{k}
\mathcal{M}_k(\partial K) \rho^k,
\end{equation}
while (\cite[(13.44)]{santalo})
\begin{equation}
\label{vrho}
\vol K\rho = \vol K + \int_0^\rho \vol_{n-1} K_\tau d \tau.
\end{equation}
\item Let $G_{n, r}$ be the Grassmanian of affine $r$-planes in
  $\mathbb{E}^n,$ with a suitably normalized invariant measure $\mu.$
Then (\cite[(14.1)]{santalo})
\begin{equation}
\label{gmeas}
\mu\left(\{L_r \in G_{n, r} \quad | \quad L_r \cap K \neq \emptyset\}\right) =
\dfrac{\omega_{n-2} \dots \omega_{n-r-1}}{(n-r)\omega_{r-1} \dots \omega_0}
\mathcal{M}_{r-1} (\partial K),
\end{equation}
where $\omega_k$ is the surface area of the unit sphere in $\mathbb{E}^{k+1}.$
\end{enumerate}
The above relationships can be used to define integral mean curvatures for
not-necessarily-$C^2$ convex bodies, and from now on the assumption of
regularity will be dropped. An important corollary of \eqref{gmeas} is
\begin{theorem}[Archimedes' axiom]
\label{arch}
$\mathcal{\partial K}$ is \emph{monotonic} under inclusion. That is, if $K_1
  \subseteq K_2,$ then $\mathcal{M}_i(\partial K_1) \leq
  \mathcal{M}_i(\partial K_2).$ The inequality is strict if $K_2 \backslash K_1$ has
  non-empty interior, and $i < n-1,$ where $n$ is the dimension of the ambient
  Euclidean space. 
\end{theorem}
\begin{remark}
I use the name ``Archimedean axiom'' because Archimedes needed a result like
this in the plane in order to make rigorous his computation of the arclengths
of curves (the circle, for example). Archimedes was unable to prove this
result from first principles so he postulated it as an axiom.
\end{remark}
\section{Polyhedra}
The formula \eqref{arho} can be used to compute integral mean curvatures for
polytopes. To wit, let $P$ be a convex polytope in $\mathbb{E}^n.$
\begin{equation}
\label{polymean}
\binom{n-1}{i}\mathcal{M}_i(\partial P) = \sum_{\mbox{codim. $i$ faces
 $f$ of $\partial P$}}
\vol_{n-i-1} f \alpha^*(f),
\end{equation}
where $\alpha^*(f)$ is the \emph{exterior angle} at $f$, described as follows:
Consider the Gauss
map, which maps each point $p$ of $\partial P$ to the set of outer normals to
the support planes to $P$ passing through $p.$  The image of all of $\partial
P$ will be the unit sphere $\mathbb{S}^{n-1} \subset \mathbb{E}^n,$ and the
combinatorial structure of $\partial P$ will induce a \emph{dual} cell
decomposition $C$ of
$\mathbb{S}^{n-1},$ in particular, the image of a codimension-$i$ face $f$ of
$\partial P$ will be an $i$-dimensional totally geodesic face $f^*$ of $C.$
The $i$-dimensional area of that face is the exterior angle at $f.$

Suppose now that translates of $P$ \emph{tile} $\mathbb{E}^n,$ so that $P$ is
a fundamental domain for a free action of a group $G$ of translations  on
$\mathbb{E^n}.$ In particular, $G$ acts on $\partial P,$ preserving the
combinatorial structure. Let the quotient by $\partial P_G.$ Then
\begin{equation}
\label{equivariant}
\binom{n-1}{i} \mathcal{M}_i = \omega_i \sum_{\mbox{codimension $i$ faces $f$ of
    $\partial P_G$}} \vol_{n-1-i} f.
\end{equation}
In particular, if $P$ is a rectangular parallelopiped:
\[
P = [0, l_1] \times [0, l_2] \times \cdots \times [0, l_n],
\]
then we obtain
\begin{equation}
\label{parallel}
\mathcal{M}_i(P) = \omega_i s_{n-1-i}(l_1, \dots, l_n).
\end{equation}
\section{Back to the Ellipsoid}
Consider again the ellipsoid
\[
E = \left\{x \quad | \quad \sum_{i=1}^n \dfrac{x_i^2}{a_i}^2 = 1\right\}.
\]
The ellipsoid $E$ is \emph{inscribed} into a \emph{box} $P$: a right
parallelopiped whose sides are parallel to the coordinate axes and have
lenghts $2a_1, \dots, 2a_n.$ Consider now the diagonal matrix $A$ whose
entries are $A_{ii} = 1/a_i.$ $A(E)$ is the unit ball $B_n,$ while $A(P)$ is
the cube $[-1, 1]^n,$ circumscribed around $B_n.$ It is clear that the cube
$[-1/\sqrt[n], 1/sqrt[n]]^n$ is inscribed in $\partial B_n,$ and so
$P/\sqrt{n}$ is inscribed in $E.$ It follows that 
\begin{equation}
\label{pinch}
\dfrac{\mathcal{M}_i(P)}{\left(\sqrt{n}\right)^{n-1-i}} =
\mathcal{M}_i\left(\dfrac{P}{\sqrt{n}}\right) \leq \mathcal{M_i}(E) \leq
\mathcal{M}_i(P),
\end{equation}
or
\begin{equation}
\label{pinch2}
\left(\dfrac{2}{\sqrt{n}}\right)^{n-1-i} \dfrac{s_{n-1-i}(l_1, \dots,
  l_n)}{\binom{n-1}{i}} \leq \mathcal{M}_i(E) \leq   2^{n-1-i}
  \dfrac{s_{n-1-i}(l_1, \dots, l_n)}{\binom{n-1}{i}}.
\end{equation}
The following theorem follows immediately.
\begin{theorem}
\label{parest}
Let 
\[f(\rho) = \frac{1}{\rho} \left[\prod_{i=1}^n(\rho + 2 l_i) - \rho^n -
  2^n \prod_{i=1}^n l_i \right].\]
Then the area of the equidistant surface $\partial E_\rho$ is bounded as follows:
\[
c_n f(\rho) \leq \vol_{n-1} \partial E_\rho \leq C_n f(\rho),
\]
where 
\[
C_n/c_n \leq \left(\sqrt{n}\right)^{n-1} \dfrac{\Gamma\left((n+1)/2\right)}{2
  \pi^{(n+1)/2}}. 
\]
\end{theorem}
\begin{remark}
\label{parestv}
A related estimate for the volume of $E_\rho$ immediately follows from Theorem
\ref{parest} and Eq. \eqref{vrho}.
\end{remark}
\section{Applications and Comments}
\subsection{Lattice Points}
The question at the root of this note was the following: Consider an ellipsoid
$E,$ and consider the number $N(e)$ of points of the integer lattice in $E.$
How do we estimate $\Delta(E) = |N(E) - \vol E|?$ 

By a standard argument (see, eg, Landau \cite{landent}), $\Delta(E)$ is
dominated by the volume of a tubular neighborhood of radius $\sqrt(n)$ of
$\partial E.$ The results of Theorem \ref{parest} and Remark \ref{parestv}, we
see that 
\begin{equation}
\label{latest}
\Delta(E) \leq C_n \int_0^{\sqrt{n}} f(\rho) d\rho \leq C_n^\prime f(\sqrt{n}),
\end{equation}
where $f$ is defined in the statement of Theorem \ref{parest}.
\subsection{Convex Bodies} It is a well-known theorem of Fritz John that for
any convex body $K,$ there exists an ellipsoid $E_K,$ such that $E_K/n
\subset K \subset E_K$ -- for centrally-symmetric $K,$ $E_K/n$ can be improved
to $E_K/\sqrt{n}.$ The results above allow us to estimate the symmetric
functions of the semi-axes (and hence the semi-axes themselves) of the John
ellipsoid in terms of the mean curvature integrals of $K$, and \emph{vice
  versa.}
\subsection{Improvements} In the paper \cite{rellip2}, among other results, the
dimensional constant for the surface area of an ellipsoid is tightened quite
considerably (from roughly $n^{n/2}$ to roughly $n.$).

\bibliographystyle{amsplain}

\end{document}